\documentclass[12pt]{amsart}
\usepackage{amssymb,amsmath,amsfonts,color,soul}
\usepackage[all]{xy}
\usepackage{hyperref}
\usepackage{upref}

\newtheorem{Theorem}{Theorem}[section] \newtheorem{Corollary}[Theorem]{Corollary} \newtheorem{Lemma}[Theorem]{Lemma} \newtheorem{Proposition}[Theorem]{Proposition} \theoremstyle{definition}
 \newtheorem{Definition}[Theorem]{Definition} \newtheorem{Example}[Theorem]{Example}  \newtheorem{Remark}[Theorem]{Remark}  \newtheorem{Question}[Theorem]{Question}    \numberwithin{equation}{section}

\newtheorem{cons}{Construction}[section]

\DeclareMathOperator{\codim}{codim}
\DeclareMathOperator{\Char}{char}

 \DeclareMathOperator{\Tor}{Tor}

 \DeclareMathOperator{\im}{im}
 \DeclareMathOperator{\rank}{rk}

\newcommand{\qis}{\simeq}
\newcommand{\del}{\partial}
\newcommand{\m}{\mathfrak{m}}
\newcommand{\ZZ}{\mathbb{Z}}

\numberwithin{equation}{section}

\begin{document}

\title[On the Rank of Multigraded Differential Modules]
{On the Rank of Multigraded\\ Differential Modules}
\author{Adam Boocher}%
\address{A.\ Boocher,
  University of San Diego, San Diego, California, USA}
\email{aboocher@sandiego.edu and aboocher@gmail.com}
\author{Justin W.~DeVries}
\address{J. DeVries, USA}
\email{s-jdevrie3@math.unl.edu}
\date{\today}

\begin{abstract}
A $\mathbb{Z}^d$-graded differential $R$-module is a
$\mathbb{Z}^d$-graded $R$-module $D$ with a morphism $\delta:D\to
D$ such that $\delta^2=0$. For $R=k[x_1,\dotsc,x_d]$, this paper
establishes a lower bound on the rank of such a differential module when
the underlying $R$-module is free. We define the
Betti number of a differential module and use it to show that when the
homology $\ker\delta/\im\delta$ of $D$ is non-zero and finite
dimensional over $k$ then there is an inequality $\rank_R D\geq 2^d$.
\end{abstract}

\maketitle

\section{Introduction}

%Let $k$ be a field and set $R=k[x_1,\dotsc,x_d]$. 
A \emph{differential
$R$-module} $D$ is an $R$-module with a square-zero
homomorphism $\delta:D\to D$ called the \emph{differential}. The homology of
$D$ is defined in the usual way: $H(D)=\ker\delta/\im\delta$.  Differential
modules have played an important role in the work of Avramov, Buchweitz,
Iyengar, and Miller on the homology of finite free complexes
\cite{MR2308849,MR2592508}. In this context, differential modules arise
naturally when working with DG-modules: some constructions with desirable
properties do not respect the grading but do preserve the differential (see
\cite{MR2592508,MR2231960} for instances of this phenomenon).

This paper is motivated by a conjecture of Avramov, Buchweitz and Iyengar, concerning differential modules over a local ring $R$ of dimension $d$.  They conjectured that if $F$ is a differential module admitting a finite free flag, such that the homology $H(F)$ has finite length then $\rank_R F\geq 2^d$ (\cite[Conjecture~5.3]{MR2308849}).

In this conjecture, a \emph{free flag} on a differential module is a
certain kind of filtration with free factors compatible with the differential (see
Definition~\ref{flag-defn}). It provides the appropriate lifting properties
for the category of differential modules. In 
\cite[Theorem~5.2]{MR2308849} the conjecture was proven when $d\leq 3$.  More generally, the conjecture is false, even for complexes!  Indeed, 
Iyengar and Walker \cite{IW} have found complexes $F$ with $H(F) \cong k\oplus k$ but 
$\rank_R F < 2^d$ for all $d \geq 8$ (provided $\Char k \neq 2$.)   However, in the 
case of resolutions, with $\Char k \neq 2$ Walker has shown \cite{W} that if $M$ is a module of finite length over $R$, then $\sum \beta_i(M) \geq 2^d$.  Further, he proves that equality holds if and only if $M$ is isomorphic to $R$ modulo a maximal $R$-sequence.

%%Phrased for a $\mathbb{Z}^d$-graded polynomial ring, this reads:
%\begin{Conjecture}\label{hor-conj}
%%Let $R=k[x_1,\dotsc,x_d]$ be a $\mathbb{Z}^d$-graded polynomial ring of dimension $d$, 
%\hl{Let $R$ be a local ring of dimension $d$ 
%and $F$ a $\mathbb{Z}^d$-graded differential $R$-module
%admitting a finite free flag. If $H(F)$ has non-zero finite length, then}
%{\[\rank_R F\geq 2^d.\]}
%\end{Conjecture}

The main result of this paper proves that the conjecture of Avramov, Buchweitz and Iyengar holds in the multigraded case when $R$ is a polynomial ring and and the differential $\delta$ has degree zero:
\begin{Theorem}\label{intro-main-thm}
Let $k$ be a field and set $R = k[x_1,\ldots, x_n]$. Let 
$F$ be a finitely generated $\mathbb{Z}^d$-graded differential
$R$-module with differential $\delta:F\to F$ that is homogeneous of degree
zero, such that $F$ is free as an $R$-module. If $H(F)$ has non-zero finite
length then $\rank_R F\geq 2^d$.
\end{Theorem}

This result is new even for complexes of $R$-modules. Given a complex of
$\mathbb{Z}^d$-graded free $R$-modules
\[F=\xymatrix{\dotso\ar[r] & F_2\ar[r] & F_1\ar[r] & F_0\ar[r] & \dotso}\]
the module $\bigoplus_i F_i$ with differential $\delta=\bigoplus_i
\partial^i$ forms a differential module. When $H(F)$ has non-zero finite length
as an $R$-module then we conclude that
\begin{equation}\label{intro-cx-conseq}
\sum_i \rank_R F_i\geq 2^d.
\end{equation}
This inequality is already known when $F$ is a resolution---i.e.\ $F_i=0$
for $i<0$ and $H_i(F)=0$ for $i\neq 0$---from the work of
Charalambous and Santoni on the Buchsbaum-Eisenbud-Horrocks problem
\cite{MR1094254,MR1028267}. Recall that for a $\mathbb{Z}$-graded
polynomial ring, the Buchsbaum-Eisenbud-Horrocks problem is to show that
$\beta_i(M)\geq\binom{d}{i}$ all $\mathbb{Z}$-graded $R$-modules $M$ with
non-zero finite length, where $\beta_i(M)$ is the $i$-th Betti number of
$M$ \cite{MR0453723,MR544153}. Summing the binomial coefficients gives
\eqref{intro-cx-conseq} when $F$ is a free resolution of a non-zero finite
length $\mathbb{Z}^d$-graded module $M$. However, when $F$ is not acyclic
it is not clear how to establish \eqref{intro-cx-conseq} without using
differential modules.

Working with differential modules provides an advantage in that it
simultaneously treats the case of free resolutions and free complexes with
homology spread among several homological degrees, as well as other
contexts. One such application arises in the conjectures of Carlsson and
Halperin concerning a lower bound on the rank of DG-modules with non-zero
finite length homology \cite{MR874170,MR787835}. For this connection between
differential modules and DG-modules see \cite[\S 5]{MR2308849}.

Some techniques available for complexes can be directly adapted to the case
of differential modules, however there are subtle difficulties that
appear. For example, there may be no way to minimize a resolution in the
category of differential modules and this creates an obstruction to applying
the usual tools of complexes (see Example~\ref{scorpion}, or
Theorem~\ref{non-pos-minimal-flag} for some positive results). Not many
techniques are available for working with differential modules. This work
should be seen as a contribution in that direction.

In \S\ref{diff-mod-section} we recall the theory of differential modules
and define a notion of a Betti number for differential modules. To place
differential modules in context, the next two sections examine connections
between differential modules and chain complexes of $R$-modules. A method
for constructing differential modules from complexes is explored in
\S\ref{compression-section}, while \S\ref{non-pos-section} presents a
result and examples to illustrate some of the difficulties in working
with differential modules as opposed to complexes. Section \ref{high-low
section} develops the main tool used for establishing lower bounds on the
rank of a differential module by adapting an inequality of Santoni
\cite{MR1028267} on the Betti numbers of a $\mathbb{Z}^d$-graded module to
a lower bound on the Betti number of $\mathbb{Z}^d$-graded differential
modules. The main result, Theorem \ref{total-Betti-bound}, is proved in
\S\ref{multi-grade-section}.  Theorem \ref{intro-main-thm} then follows. 

\subsection{Related Results} If $M$ is isomorphic to $R$ modulo a maximal $R$-sequence 
then $M$ is resolved by the Koszul complex, and thus $\beta_i(M) = {d \choose i}$.
In the $\mathbb{Z}^d$-graded setting, for resolutions, Charalambous and Evans \cite{C,CE} showed 
that if $M$ is a multigraded module of finite length and $M$ is not isomorphic 
to $R$ modulo a maximal $R$-sequence then at least one of the following holds:
$$(a) \ \ \ \ \  \beta_i(M) \geq {d \choose i} + {d -1 \choose i-1}, \mbox{ for all $i$},$$
$$(b) \ \ \ \ \  \beta_i(M) \geq {d \choose i} + {d -1 \choose i}, \mbox{ for all $i$}.$$
These inequalities are false if one drops the multigraded hypothesis.  (Consider $M = k[x,y,z]/(x^2-y^2, y^2-z^2, xy,yz)$.)  If $M$ is not finite length, 
then one can reformulate bounds on Betti numbers by replacing $d$ with $c = \codim M$.  However, even if $M = S/I$ is a cyclic module, neither
(a) nor (b) need hold.  (Consider $M = k[x,y,z,u,v]/(xy,yz,zu,uv,vx)$.)  Nevertheless, in \cite{BS}, it was shown that $\sum \beta_i(S/I) \geq 2^c + 2^{c-1}$ for any monomial ideal of height $c$ that does not define a complete intersection.   For related results, see also \cite{CEM,BoocherW}. 
These results suggest the following:
\begin{Question}\ 
\begin{itemize}
\item Suppose $M$ is a $\ZZ^d$-graded $R$-module ($M$ not necessarily of finite length) that is not isomorphic to $R$ modulo a maximal $R$-sequence. Is $\sum \beta_i(M) \geq 2^c + 2^{c-1}$ where $c = \codim M$?
\item Let $F$ be a finitely generated $\mathbb{Z}^d$-graded differential
$R$-module with differential $\delta:F\to F$ that is homogeneous of degree
zero, such that $F$ is free as an $R$-module. If $H(F)$ has non-zero finite
length with $H(F)$ not isomorphic to $R$ modulo a regular sequence then is $\rank_R F \geq 2^d + 2^{d-1}?$
\end{itemize}
\end{Question}

\section{Differential modules}\label{diff-mod-section}

Throughout, $k$ is a field, $R=k[x_1,\dotsc,x_d]$ is the
standard $\mathbb{Z}^d$-graded polynomial ring and $\mathfrak{m}=(x_1,\dotsc,x_d)$.
To be specific, the grading on $R$ is such that the degree
$\deg(x_i)\in\mathbb{Z}^d$ of variable $x_i$ is $(0,\dotsc,0,1,0,\dotsc,0)$
with the 1 appearing in the $i$-th coordinate. For
$\mathbf{m}\in\mathbb{Z}^d$, we write $\mathbf{m}_i$ to denote the $i$-th
coordinate. Two elements $\mathbf{a},\mathbf{b}\in\mathbb{Z}^d$ are
compared coordinate-wise by setting $\mathbf{a}\leq\mathbf{b}$ if
$\mathbf{a}_i\leq\mathbf{b}_i$ for all $i$. This makes $\mathbb{Z}^d$ into
a partially ordered group.

Recall that a $\mathbb{Z}^d$-graded module $M$ over $R$ is an $R$-module
that has a decomposition $\bigoplus_{\mathbf{m}\in\mathbb{Z}^d}
M_\mathbf{m}$ as abelian groups such that multiplication by an element of
$R$ of degree $\mathbf{n}$ takes $M_\mathbf{m}$ to
$M_{\mathbf{m}+\mathbf{n}}$. An $R$-linear map $\phi$ between
$\mathbb{Z}^d$-graded modules $M$ and $N$ is a
\emph{morphism} if $\phi(M_\mathbf{m})\subseteq N_\mathbf{m}$. In
particular, a complex of $\mathbb{Z}^d$-graded modules is required to have
morphisms for its differentials.

For $\mathbf{d}\in\mathbb{Z}^d$ the shifted (or twisted) module
$M(\mathbf{d})$ is defined to be $M_{\mathbf{d}+\mathbf{m}}$ in degree
$\mathbf{m}$ for each $\mathbf{m}\in\mathbb{Z}^d$, with the same $R$-module
structure as $M$. Given a morphism $\phi:M\to N$ the shifted morphism
$M(\mathbf{d})\to N(\mathbf{d})$ defined by $x\mapsto \phi(x)$ is denoted
$\phi(\mathbf{d})$.

We will work with $\mathbb{Z}^d$-graded modules and
$\mathbb{Z}^d$-graded differential modules, so definitions will be given
in that context for simplicity; see \cite{MR2308849,MR1731415,my-thesis}
for details concerning arbitrary differential modules.
%add something here about cited references working for multigraded?

\begin{Definition}\label{diff-mod-defn}
A \emph{$\mathbb{Z}^d$-graded differential $R$-module with differential degree
$\mathbf{d}\in\mathbb{Z}^d$} is a $\mathbb{Z}^d$-graded $R$-module $D$
with a morphism $\delta:D\to D(\mathbf{d})$ such that the composition
\[\xymatrix{D(-\mathbf{d})\ar^-{\delta(-\mathbf{d})}[r] &
D\ar^-{\delta}[r] & D(\mathbf{d})}\] is zero. We say that $\delta$ is the
\emph{differential} of $D$.

When $D$ and $E$ are $\mathbb{Z}^d$-graded differential modules with the
same differential degree, define a \emph{morphism} $\phi:D\to E$
to be a morphism of $\mathbb{Z}^d$-graded modules satisfying
$\delta^E\circ\phi=\phi\circ\delta^D$. For a fixed differential degree, the
category of $\mathbb{Z}^d$-graded differential modules with this notion of
a morphism is an abelian category.

The homology of a differential module $D$ is the $\mathbb{Z}^d$-graded
$R$-module
\[H(D)=\ker\delta/\im(\delta(-\mathbf{d})).\]
The $\mathbb{Z}^d$-grading on $H(D)$ is inherited from $D$ by considering
$\ker\delta$ and $\im(\delta(-\mathbf{d}))$ as submodules of $D$ with the
induced grading. Any $\mathbb{Z}^d$-graded $R$-module, in particular
$H(D)$, will be considered as a differential module with zero differential.

In the usual way, a morphism $\phi:D\to E$ induces a map in homology
$H(\phi):H(D)\to H(E)$. If $H(\phi)$ is an isomorphism we say that $\phi$
is a \emph{quasi-isomorphism} and write $D\qis E$ or
$\phi:D\xrightarrow{\qis} E$. Given an exact sequence of differential
modules
\[\xymatrix{0\ar[r] & D_1\ar^{\alpha}[r] & D_2\ar^{\beta}[r] & D_3\ar[r] &
0}\]
there is an induced long exact sequence of in homology,
\[
\xymatrix{\dotso\ar[r]
& H(D_1)(i\mathbf{d})\ar^{H(\alpha)(i\mathbf{d})}[r]
& H(D_2)(i\mathbf{d})\ar^{H(\beta)(i\mathbf{d})}[r]
& H(D_3)(i\mathbf{d})\ar^-{\gamma(i\mathbf{d})}[r]
& H(D_1)((i+1)\mathbf{d})\ar[r] & \dotso}
\]
where $i$ ranges over the integers, and each map is a morphism of
$\mathbb{Z}^d$-graded modules (in particular, has degree $\mathbf{0}$).
We summarize this sequence by the following diagram
\begin{equation}\label{long-exact-triangle}
\xymatrix{H(D_1)\ar^{H(\alpha)}[rr] && H(D_2)\ar^{H(\beta)}[dl]\\
& H(D_3)\ar^*!/^3pt/{\labelstyle \gamma}[ul]|{\bigcirc}\\}
\end{equation}
where the circle indicates that $\gamma$ is a homomorphism of degree
$\mathbf{d}$.

See \cite[Chap.~IV \S 1]{MR1731415} for a proof.
\end{Definition}

Bounds on the rank of a differential module will be obtained by comparing
the rank and an invariant that we call the Betti number of a differential
module. To define the Betti number we will need a notion of a tensor
product of differential modules. However, adapting the usual definition of
a tensor product between complexes fails to produce a differential module
when applied to two differential modules. To work around this we recall the
construction of a tensor product of a complex and a differential module,
along with some of its properties \cite[\S 1]{MR2308849}.

\begin{Definition}
For a complex $C$ of $\mathbb{Z}^d$-graded $R$-modules and a
$\mathbb{Z}^d$-graded differential
$R$-module $D$ with differential degree $\mathbf{d}$, define a
$\mathbb{Z}^d$-graded differential module $C\boxtimes_R D$ by setting
\[C\boxtimes_R D=\bigoplus_{i\in\mathbb{Z}} (C_i(-i\mathbf{d})\otimes_R
D),\]
with differential defined by
\[\delta^{C\boxtimes_R D}(c\otimes d)=\del^C(c)\otimes d + (-1)^i
c\otimes\delta^D(d),\]
for $c\otimes d\in C_i(-i\mathbf{d})\otimes_R D$. This makes $C\boxtimes_R
D$ into a $\mathbb{Z}^d$-graded differential $R$-module with differential
degree $\mathbf{d}$.
\end{Definition}

We will need the following facts concerning this product. These results are
proved in \cite{MR2308849} for arbitrary differential modules, but the
proofs hold for $\mathbb{Z}^d$-graded differential modules with the obvious
modifications.

\begin{Proposition}[{\cite[1.9.3]{MR2308849}}]\label{box-tensor-assoc}
Let $X$ and $Y$ be $\mathbb{Z}^d$-graded complexes and let $D$ be a
$\mathbb{Z}^d$-graded differential module. Then there is a natural
isomorphism of $\mathbb{Z}^d$-graded differential modules:
\[(X\otimes_R Y)\boxtimes_R D=X\boxtimes_R(Y\boxtimes_R D).\]
\end{Proposition}

\begin{Proposition}[{\cite[Proposition~1.10]{MR2308849}}]\label{cx-tensor-exact}
Let $X$ and $Y$ be bounded below $\mathbb{Z}^d$-graded complexes of flat
$R$-modules, i.e.\ $X_i=Y_i=0$ for sufficiently small $i$. Then
\begin{enumerate}
\item\label{enum:cx-tensor-exact-1}
the functor $X\boxtimes_R -$ preserves exact sequences and
quasi-isomorphisms;
\item\label{enum:cx-tensor-exact-2}
a quasi-isomorphism $\phi:X\to Y$ induces a quasi-isomorphism
\[\phi\boxtimes_R D: X\boxtimes_R D\to Y\boxtimes_R D\]
for all $\mathbb{Z}^d$-graded differential $R$-modules $D$.
\end{enumerate}
\end{Proposition}

Using this tensor product, we can define a $\Tor$ functor between
$R$-modules and differential $R$-modules, and hence define a Betti number.

\begin{Definition}
For a $\mathbb{Z}^d$-graded differential $R$-module $D$ and a
$\mathbb{Z}^d$-graded $R$-module $M$ set
\[\Tor^R(M,D)=H(P\boxtimes_R D)\]
where $P$ is a $\mathbb{Z}^d$-graded free resolution of $M$. This is
well-defined as different choices of free resolution produce
quasi-isomorphic differential modules by Proposition~\ref{cx-tensor-exact}.
\end{Definition}

\begin{Definition}\label{Betti-defn}
We define $\beta^R_\mathbf{m}(D)$ to be the \emph{Betti number in degree
$\mathbf{m}\in\mathbb{Z}^d$} of a differential $R$-module $D$:
\[\beta^R_\mathbf{m}(D)=\rank_k\Tor^R(k,D)_\mathbf{m}.\]
Summing over all degrees gives the \emph{Betti number} $\beta^R(D)$:
\[\beta^R(D)=\sum_{\mathbf{m}\in\mathbb{Z}^d} \beta^R_\mathbf{m}(D)
=\rank_k\Tor^R(k,D).\]
\end{Definition}

The connection between ranks of differential modules and Betti numbers is
provided by free flags, a notion of a free resolution for differential
modules \cite[\S 2]{MR2308849}.

\begin{Definition}\label{flag-defn}
A \emph{free flag} on a differential module $F$ is a family
$\{F^n\}_{n\in\mathbb{Z}}$ of $\mathbb{Z}^d$-graded $R$-submodules such
that
\begin{enumerate}
\item $F^{n}=0$ for $n<0$,
\item $F^n\subseteq F^{n+1}$ for all $n$,
\item $\delta^F(F^{n+1})\subseteq F^{n}$ for all $n$,
\item $\bigcup_{n\in\mathbb{Z}} F^n=F$,
\item $F^{n+1}/F^{n}$ is a free $R$-module for all $n$.
\end{enumerate}

A $\mathbb{Z}^d$-graded differential module $F$ with a free flag
\emph{resolves} $D$ if there is a quasi-isomorphism
$\xymatrix{F\ar^{\qis}[r] & D}$ in the category of $\mathbb{Z}^d$-graded
differential modules.
%do we want to require surjectivity?
\end{Definition}

\begin{Remark}\label{345free-rmk}
Properties (3), (4), (5) in Definition {\ref{flag-defn}} imply that for each $n$, $F^{n+1} = F^{n} \oplus F_{n+1}$ where $F_{n+1}$ is a free module with $\delta^F (F_{n+1}) \subseteq F^{n}$. 
\end{Remark}

Many properties of free bounded-below complexes have analogs for differential modules with
free flags. We will use the following two.

\begin{Proposition}\label{free-flag-lifting}
Let $D_1$ and $D_2$ be differential modules and let $F$ be a differential
module with a free flag. If $\alpha:D_1\to D_2$ is a surjective
quasi-isomorphism and $\beta:F\to D_2$ is a morphism then there is a
morphism $\gamma:F\to D_1$ such that the following diagram commutes:
\[\xymatrix{ & D_1\ar_{\qis}^{\alpha}[d] \\ F\ar_{\beta}[r]
\ar@{-->}^{\gamma}[ur] & D_2}\]
\end{Proposition}

\begin{proof}[Sketch of proof]
Let $\{F^n\}_{n\in\mathbb{Z}}$ be a free flag on $F$. Define $\gamma:F\to
D_1$ inductively by defining $\gamma^n:F^n\to D_1$. We can define
$\gamma^0:F^0\to D_1$ using the usual lifting properties since $F^0$ is a
free $R$-module. For $n>0$, by Remark {\ref{345free-rmk}} we have $F^n=F^{n-1}\oplus F_n$ for a free module $F_n$ with $\delta^F(F_n) \subseteq F^{n-1}$.  Assuming that we have $\gamma^{n-1}:F^{n-1}\to D_1$
defined, we can define $\gamma^n:F^n\to D_1$ by using the lifting
properties of the free module $F_n$ to define a map $F_n\to
D_1$. The lifting used is important since we need
$\delta^{D_1}\gamma=\gamma\delta^F$ and $\alpha\gamma=\beta$. However, any
lifting can be modified by adding an appropriate boundary of $D_1$ so that
it has the desired properties.
\end{proof}

\begin{Proposition}[{\cite[Proposition~2.4]{MR2308849}}]\label{diff-mod-tensor-exact}
Let $F$ be a $\mathbb{Z}^d$-graded differential module with a free flag.
Then the functor $-\boxtimes_R F$ preserves exact sequences and
quasi-isomorphisms.
\end{Proposition}

With differential modules that admit a free flag providing a resolution
of a differential module, the $\Tor$ functor is balanced, which gives the
connection between the rank and Betti number of a differential module.

\begin{Lemma}\label{tor-balanced}
Let $P$ be a free resolution of a $\mathbb{Z}^d$-graded module $M$ and
let $F$ be a free flag resolving a $\mathbb{Z}^d$-graded differential
module $D$. Then $H(P\boxtimes_R D)$ is isomorphic to $H(M\boxtimes_R F)$
as $\mathbb{Z}^d$-graded $R$-modules.
\end{Lemma}

\begin{proof}
Let $\varepsilon:P\to M$ and $\eta:F\to D$ be $\mathbb{Z}^d$-graded
quasi-isomorphisms. Then there are $\mathbb{Z}^d$-graded morphisms
\[\xymatrix{P\boxtimes_R D & P\boxtimes_R F\ar_{P\boxtimes_R \eta}[l]
\ar^{\varepsilon\boxtimes_R F}[r] & M\boxtimes_R F.}\]
By Proposition~\ref{cx-tensor-exact} and
Proposition~\ref{diff-mod-tensor-exact} these are quasi-isomorphisms.
\end{proof}

In particular $\beta_m^R(F) = \rank_k H(k\otimes_R F)_m.$

\begin{Theorem}\label{Betti-rank-thm}
Let $F$ be a $\mathbb{Z}^d$-graded differential module admitting a free flag.
For all degrees $\mathbf{m}\in\mathbb{Z}^d$ we have
$\beta^R_\mathbf{m}(F)\leq\rank_k F_\mathbf{m}$. Therefore,
\[\beta^R(F)\leq\rank_R F.\]
\end{Theorem}

\begin{proof}
By Lemma~\ref{tor-balanced},
\[\beta^R_\mathbf{m}(F)=\rank_k \Tor^R(k,F)_\mathbf{m}
=\rank_k H(k\boxtimes_R F)_\mathbf{m}.\]
Since $k$ is an $R$-module, $k\boxtimes_R F=k\otimes_R F$. Since
$H(k\otimes_R F)_\mathbf{m}$ is a subquotient of $(k\otimes_R
F)_\mathbf{m}$, we have
\[\beta^R_\mathbf{m}(F)=\rank_k H(k\boxtimes_R F)_\mathbf{m}\leq \rank_k
(k\otimes_R F)_\mathbf{m}.\]
Summing over all degrees gives the inequality for the Betti number,
\[\beta^R(F)=\sum_{\mathbf{m}\in\mathbb{Z}^d} \beta^R_\mathbf{m}(F)
\leq \sum_{\mathbf{m}\in\mathbb{Z}^d} \rank_k (k\otimes_R F)_\mathbf{m}
=\rank_k k\otimes_R F=\rank_R F.\qedhere\]
\end{proof}

\begin{Remark}\label{Betti-minl-rmk}
When $\delta(F)\subseteq\m F$ we have
$\beta^R_\mathbf{m}(F)=\rank_k F_\mathbf{m}$ as the differential of
$k\boxtimes_R F$ is zero. In general, the inequality can be strict; see
Example~\ref{scorpion}.
\end{Remark}

We finish this section by recording a property of the $\Tor$ functor for
use later.

\begin{Lemma}
Consider an exact sequence of $\mathbb{Z}^d$-graded differential
$R$-modules
\[\xymatrix{0\ar[r] & D_1\ar^{\alpha}[r] & D_2\ar^{\beta}[r] &  D_3\ar[r] &
0}.\]
For each $\mathbb{Z}^d$-graded $R$-module $M$ there is an exact commutative
diagram of $\mathbb{Z}^d$-graded differential modules:
\[\xymatrix{ \Tor^R(M,D_1) \ar^-{\Tor(M,\alpha)}[rr] &&
\Tor^R(M,D_2)\ar^-*!/^2pt/{\labelstyle \Tor(M,\beta)}[dl]\\
& \Tor^R(M,D_3)\ar^*!/^3pt/{\labelstyle \gamma}[ul]|{\bigcirc}\\}\]
\end{Lemma}

\begin{proof}
Take a free resolution $P$ of the module $M$. By
Proposition~\ref{cx-tensor-exact} the sequence of differential modules
remains exact after applying $P\boxtimes_R -$:
\[\xymatrix{0\ar[r] & P\boxtimes_R D_1\ar^{P\boxtimes\alpha}[r] &
P\boxtimes_R D_2\ar^{P\boxtimes\beta}[r] & P\boxtimes_R D_3\ar[r] & 0}.\]
The diagram \eqref{long-exact-triangle} coming from this exact
sequence is the desired one.
\end{proof}

\section{Compression}\label{compression-section}

Every complex of $R$-modules produces a differential module by forming its
\emph{compression}. This construction allows results about differential
modules to be translated to results about complexes of modules. In fact,
the differential modules produced by compressing always have differential
degree $\mathbf{0}$ so it is sufficient to restrict to differential modules
with differential degree $\mathbf{0}$ if one is interested in
establishing results about complexes. Note that not every differential
module of differential degree $\mathbf{0}$ arises this way (see
Example~\ref{deg-0-scorpion}).

\begin{cons}[{\cite[1.3]{MR2308849}}]
If $C$ is a complex of $\mathbb{Z}^d$-graded $R$-modules, then its
\emph{compression} is the $\mathbb{Z}^d$-graded differential module
\[C_\Delta=\bigoplus_{i\in\mathbb{Z}} C_i\]
with differential $\delta^{C_\Delta}=\bigoplus_{i\in\mathbb{Z}} \del^C_i$.

We have $\deg(\delta^{C_\Delta})=\mathbf{0}$ because the differentials of
the complex $C$ are required to have degree zero. By the definition of
$\delta^{C_\Delta}$, we have $H(C_\Delta)=\bigoplus_{i\in\mathbb{Z}}
H_i(C)$.

When the complex $C$ is bounded below and consists of free $R$-modules then
the compression has a free flag. Indeed, suppose $C_i=0$ for $i$
sufficiently small. Then after an appropriate shifting, setting $F^n=\bigoplus_{i\leq n} C_i$ forms a free
flag.
\end{cons}

Computing the Betti number of a compression of a minimal complex is a straight-forward application
of Theorem~\ref{Betti-rank-thm} and Remark~\ref{Betti-minl-rmk}.

\begin{Lemma}
Let $C$ be a bounded below complex of free modules that is minimal in the
sense that $\del_n^C(C_n)\subseteq \m C_{n-1}$. Then
\[\beta(C_\Delta)=\sum_i \rank_R C_i.\]
When $C$ is a minimal free resolution of a module $M$ we have
\[\beta(C_\Delta)=\sum_i \beta_i(M),\]
where $\beta_i(M)$ is the usual Betti number of $M$.
\end{Lemma}

\begin{proof}
Since $C$ is a bounded below complex of free modules, $C_\Delta$ has a free
flag. We have
\[\delta(C_\Delta)=\bigoplus_{i\in\mathbb{Z}} \del_i(C_i)\subseteq
\bigoplus_{i\in\mathbb{Z}} \m C_{i-1}=\m C_\Delta,\]
so by Remark~\ref{Betti-minl-rmk} we have
\[\beta(C_\Delta)=\rank_R C_\Delta=\sum_i \rank_R C_i.\]
When $C$ is a minimal free resolution of $M$ we have $\rank_R
C_i=\beta_i(M)$, which completes the proof.
\end{proof}

Obviously differential modules with non-zero differential degree do not
come from compressing a complex, but the following shows that there
are also differential modules with differential degree zero that are not
compressions of a complex.

\begin{Example}\label{deg-0-scorpion}
Let $R=k[x,y]$ and let $F=R(0,0)\oplus R(-1,0)\oplus R(0,-1)\oplus
R(-1,-1)$. Viewing $F$ as column vectors, define a differential $\delta$ by
left-multiplication by the matrix
\[\begin{bmatrix} 0 & x & y & xy\\ 0 & 0 & 0 & -y\\ 0 & 0 & 0 & x\\
0 & 0 & 0 & 0\end{bmatrix}.\]
This is a differential module with $\deg\delta=\mathbf{0}$. Represented
diagrammatically this has the form of a Koszul complex on $x,y$ modified by
adding an additional map:
\begin{equation*}
\xymatrix{R(-1,-1)\ar_-{\bigl[\begin{smallmatrix} -y\\
x\\\end{smallmatrix}\bigr]}[r]
\ar@(ur,ul)^{xy}[rr] &
R(-1,0)\oplus R(0,-1)\ar_-{[\begin{smallmatrix} x & y\end{smallmatrix}]}[r] &
R(0,0)\ar[r] & 0}.
\end{equation*}
Reading the diagram from right to left produces a free flag:
\begin{multline*}
0\subset R(0,0)\subset R(0,0)\oplus R(-1,0)\oplus R(0,-1)\subset\\
R(0,0)\oplus R(-1,0)\oplus R(0,-1)\oplus R(-1,-1)=F.
\end{multline*}
To calculate $H(F)$, consider the first differential submodule of the flag
$F^0=R(0,0)$. It is straight-forward to see that
\begin{align*}
H(F^0)&=R(0,0)\\
H(F/F^0)&=(R(-1,0)\oplus R(0,-1))/R(-y\oplus x).
\end{align*}
From the short exact sequence
\[\xymatrix{0\ar[r] & F^0\ar[r] & F\ar[r] & F/F^0\ar[r] & 0}\]
we have the long exact sequence
\[\xymatrix{\dotso\ar[r] & H(F/F^0)\ar^-{\alpha}[r] &
H(F^0)\ar^-{\beta}[r] & H(F)\ar[r] &
H(F/F^0)\ar^-{\alpha}[r] & H(F^0)\ar[r] & \dotso}\]
where the map $\alpha$ is given by the matrix $\begin{bmatrix} x &
y\end{bmatrix}$. Since $\alpha$ is injective, $\beta$ must be a surjection,
giving
\[H(F)=H(F^0)/\im\alpha = R/(x,y)=k.\]

To compute the Betti number, note that $\delta(F)\subseteq\m F$, so we have
$\beta^R(F)=\rank_R F=4$ by Remark~\ref{Betti-minl-rmk}.
\end{Example}

\section{Non-positive differential degree}\label{non-pos-section}

Every differential $R$-module with a free flag is free as an
$R$-module, but not conversely (see Example \ref{minimized-scorpion}). Even
when a differential module admits a free flag there may be no way to
``minimize,'' unlike finite free complexes that can be decomposed into an
acyclic complex and a minimal complex $C$ with $\del(C)\subseteq \m C$ (see
Example \ref{scorpion}). Restricting to the case of a finitely generated differential module
$D$ with $\deg\delta^D\leq\mathbf{0}$ we can avoid both of these
difficulties.

\begin{Theorem}\label{non-pos-minimal-flag}
Let $F$ be a finitely generated $\mathbb{Z}^d$-graded differential
$R$-module with $\deg\delta^F\leq\mathbf{0}$ that is free as an $R$-module.
Then $F$ has a free flag and a submodule $F'$ that is a direct summand in
the category of $\mathbb{Z}^d$-graded differential $R$-modules such that
\begin{enumerate}
\item $F'$ has a free flag,
\item $\delta(F')\subseteq \m F'$,
\item $H(F')=H(F)$.
\end{enumerate}
\end{Theorem}

\begin{Remark}
The hypothesis that $\deg\delta\leq\mathbf{0}$ is necessary. See Examples
\ref{scorpion} and \ref{minimized-scorpion}.
\end{Remark}

\begin{proof}
We induce on $\rank_R F$: if $\rank_R F=1$ then the differential of $F$ is
multiplication by an element of $R$. Since $R$ is a domain, this element
must be zero; hence $F^0=F$ is a free flag. As $\delta(F^0)=0$ we conclude
that $\delta(F)\subseteq \m F$ as well.

Now suppose $\rank_R F>1$. If $\delta(F)\not\subseteq \m F$ then there is
some homogeneous basis element $e$ with $\delta(e)\not\in \m F$. We first
show that $\overline{e},\delta(\overline{e})\in F/\m F$ are linearly
independent over $k$. Suppose that there is a linear relation
$\delta(\overline{e})=a\overline{e}$ with $a\in k$ nonzero. Since $\delta^2=0$, we
have $0=a\delta(\overline{e})=a^2\overline{e}$, a contradiction.

So $\overline{e}$ and $\overline{\delta(e)}$ are linearly independent. By Nakayama's lemma
we can take $\{e,\delta(e)\}$ to be part of a basis of $F$. Let
$E=Re\oplus R\delta(e)$. Then $E$ is a differential submodule. So we
have an exact sequence of differential modules:
\begin{equation}\label{eqn:ses}
\xymatrix{0\ar[r] & E\ar[r] & F\ar[r] & F/E\ar[r] & 0}.
\end{equation}
Since $H(E)=0$, the long exact sequence in homology coming from
\eqref{eqn:ses} shows that $H(F/E)=H(F)$. The module $F/E$ is free
since $E$ is generated by basis elements of $F$.  By a standard splitting argument, and Proposition \ref{free-flag-lifting} there is a differential $R$-module $G\cong F/E$ so that $F = E \oplus G$.  So by
induction $G$ has a free flag $\{G^n\}_{n\in\mathbb{Z}}$.  Proposition \ref{free-flag-lifting} shows that 

Setting
\begin{align*}
F^0&=R\delta^F(e),\\
F^1&=R\delta^F(e)\oplus Re,\\
F^n&=R\delta^F(e)\oplus Re\oplus G^{n-2},\quad n\geq 2
\end{align*}
gives a free flag on $F$. The induction hypothesis also
shows that $G_0$ has a direct summand $F'$ with a free flag such that
$\delta(F')\subseteq \m F'$ and such that $H(F')=H(G_0)=H(F)$.  This
completes the proof when $\delta(F)\not\subseteq \m F$.

Now suppose that $\delta(F)\subseteq \m F$. In this case it suffices to
show that $F$ has a free flag. Let $e_1,\dotsc,e_n$ be a homogeneous basis
for $F$. Let $\mathbf{n}$ be a minimal element of
$\{\deg(e_1),\dotsc,\deg(e_n)\}$ under the partial order on $\mathbb{Z}^d$.
Set
\[G=\bigoplus_{\deg(e_i)=\mathbf{n}} Re_i.\]
Then $\delta^F(G)\subseteq G$ since $\deg(\delta^F(e_i))\leq\deg(e_i)$ for
all $i$ as the degree of $\delta^F$ is non-positive in each coordinate. So
$G$ is a differential submodule.

We claim that $\delta^F|_G=0$. When $\deg\delta^F<\mathbf{0}$, we have
$\delta^F|_G=0$ as $\deg(\delta^F(e_i))<\deg(e_i)$ and all the generators
$e_i$ of $G$ have the same degree. When $\deg\delta^F=\mathbf{0}$ the
matrix representing $\delta^F|_G$ has entries in $k$ since all generators
of $G$ are in the same degree. So $\delta^F|_G=0$, otherwise there would be
an element of $\delta^F(G)$ that is not in $\m F$, contrary to assumption.

Since $\delta^F|_G=0$ we get $\delta^F(F^0)=0$ by setting $F^0=G$. As $F^0$
is generated by basis elements of $F$, the quotient $F/F^0$ is a free
$R$-module and using a standard splitting argument, there is a free $R-$module $F_0 \cong F/F^0$ such that 
$F = F^0 \oplus F_0$. The induction hypothesis produces a free flag
$\{G^n\}_{n\in\mathbb{Z}}$ for $F_0$. Setting $F^n=F^0\oplus G^{n-1}$ for
$n\geq 0$ and $F^n=0$ for $n<0$ gives a free flag on $F$.
\end{proof}

The next example illustrates several difficulties in dealing with
differential modules with non-zero differential degree. It provides an
obstruction to extending Theorem~\ref{total-Betti-bound} to differential
modules with $\deg\delta>\mathbf{0}$. Furthermore, by
\cite[Theorem~5.2]{MR2308849}, a differential module over $k[x,y]$ (a ring of dimension $d=2$)
with non-zero finite length homology and a finite free flag must have rank
at least $2^d = 4$.   So this example also shows that
Theorem~\ref{non-pos-minimal-flag} cannot be extended to differential
modules with $\deg\delta>\mathbf{0}$ as no summand can have a free flag.

\begin{Example}\label{scorpion}
Let $R=k[x,y]$ and let $F=R(0,0)\oplus R(0,1)\oplus R(1,0)\oplus R(1,1)$
have differential given by the matrix,
\[\delta=
\begin{bmatrix}
0 & x & y & 1\\
0 & 0 & 0 & -y\\
0 & 0 & 0 & x\\
0 & 0 & 0 & 0\\
\end{bmatrix}.\]
This is a differential module with differential degree $(1,1)$.
As a diagram it is
\begin{equation}\label{scorpion-diagram}
\xymatrix{R(1,1)\ar_-{\bigl[\begin{smallmatrix} -y\\
x\\\end{smallmatrix}\bigr]}[r]
\ar@(ur,ul)^{1}[rr] &
R(0,1)\oplus R(1,0)\ar_-{[\begin{smallmatrix} x & y\end{smallmatrix}]}[r] &
R(0,0)\ar[r] & 0}.
\end{equation}
As in Example~\ref{deg-0-scorpion}, reading the diagram from right to left
gives a free flag. The same computation from Example~\ref{deg-0-scorpion}
shows that $H(F)=k$. As $F$ has a free flag, we can compute $\beta^R(F)$ as
$\rank_k H(k\otimes_R F)$ by the statement immediately after Lemma \ref{tor-balanced}. Applying $k\boxtimes_R -$ to
\eqref{scorpion-diagram} we have the vector space $k^4$ (suppressing the
grading) with differential given by the diagram:
\[\xymatrix{k\ar_{0}[r]\ar@(ur,ul)^{1}[rr] & k^2\ar_{0}[r] & k\ar[r] &
0}.\]
The homology is $k^2$, so $\beta^R(F)=2$.
\end{Example}

This final example shows that a differential module that is free as an
$R$-module need not have a free flag; thus Theorem~\ref{non-pos-minimal-flag}
cannot be strengthened to apply to differential modules with
$\deg\delta>\mathbf{0}$.

\begin{Example}[{\cite[Example~5.6]{MR2308849}}]\label{minimized-scorpion}
Let $F$ be as in Example~\ref{scorpion}. Let $e$ be the basis element in
degree $(-1,-1)$ and set $G=Re\oplus R\delta^F(e)$. Then a calculation shows
that $F/G$ is the differential module
$D=R(0,1)\oplus R(1,0)$ with
\[\delta=\begin{bmatrix} xy & -y^2\\ x^2 & -xy\end{bmatrix}.\]
This is a differential module with $\deg\delta=(1,1)$. Since $H(G)=0$, an
exact sequence argument shows that the map $F\to F/G$ is a
quasi-isomorphism; hence $H(D)=H(F)=k$. As $F$ admits a free flag, it is a
resolution of $D$. So we have $\beta^R(D)=\beta^R(F)=2$.

The differential module $D$ itself cannot have a free flag since $\rank_R
D=2<4$, as noted before Example~\ref{scorpion}.
\end{Example}

\section{High-low decompositions}\label{high-low section}

The main tool, Theorem~\ref{tor-thm}, we use for finding a bound on the
Betti number comes from an inequality of Santoni \cite{MR1028267}
reformulated to apply to differential modules. The essential idea is
to use information about the ``top'' and ``bottom'' degree parts to derive
information about the entire module.  The meaning of ``top'' and ``bottom''
is made precise by a \emph{high-low decomposition}; see
Definition~\ref{single-var-high-low}.

Let $y$ be an indeterminate over $R=k[x_1,\dotsc,x_d]$ with $\deg
y=(0,\dotsc,0,1)\in\mathbb{Z}^{d+1}$, so that $R[y]$ is a
$\mathbb{Z}^{d+1}$-graded ring. In this section we will be concerned with
comparing $\mathbb{Z}^{d+1}$-graded differential modules over $R[y]$ with
$\mathbb{Z}^d$-graded differential modules over $R$. Via the inclusion
$R\hookrightarrow R[y]$, any $\mathbb{Z}^{d+1}$-graded differential module
over $R[y]$ can be considered as a $\mathbb{Z}^d$-graded differential
module over $R$, with the action of $R$ fixing the $(d+1)$-th coordinate of
the $\mathbb{Z}^{d+1}$-grading. The following result allows this change of
rings to be applied to the $\Tor$ functor.

\begin{Lemma}\label{change-of-rings}
Let $M$ be a $\mathbb{Z}^{d+1}$-graded $R[y]$-module and $D$ a
$\mathbb{Z}^{d}$-graded differential $R$-module. View
$R[y]\boxtimes_R D$ as a $R[y]$-module via the action $r(s\otimes
d)=(rs)\otimes d$. Then
\[\Tor^{R[y]}(M,R[y]\boxtimes_R D)\cong \Tor^R(M,D)\]
as $\mathbb{Z}^{d+1}$-graded differential modules.
\end{Lemma}

\begin{proof}
Let $P$ be a $\mathbb{Z}^{d+1}$-graded free resolution of $M$ over $R[y]$.
Then by using Proposition~\ref{box-tensor-assoc} one gets:
\begin{align*}
\Tor^{R[y]}(M,R[y]\boxtimes_R D)&=H(P\boxtimes_{R[y]} (R[y]\boxtimes_R D))\\
&\cong H((P\otimes_{R[y]} R[y])\boxtimes_R D)\\
&\cong H(P\boxtimes_R D)\\
&=\Tor^R(M,D).\qedhere
\end{align*}
\end{proof}

Let $\mathcal{C}$ be a class of $\mathbb{Z}^{d+1}$-graded differential
$R[y]$-modules which is closed under taking submodules and quotients. Take
$\lambda$ to be a superadditive function from $\mathcal{C}$ to an ordered
commutative monoid such that $\lambda(C)\geq 0$ for all $C\in\mathcal{C}$.
Recall that $\lambda$ is superadditive if an exact sequence
\[\xymatrix{0\ar[r] & A\ar[r] & B\ar[r] & C\ar[r] & 0}\]
of differential modules in $\mathcal{C}$ gives an inequality
$\lambda(B)\geq \lambda(A)+\lambda(C)$.

\begin{Example}
For our purposes, $\mathcal{C}$ will be the collection
of finitely generated $\mathbb{Z}^{d+1}$-graded differential $R[y]$-modules with non-zero
homology in at most finitely many degrees and $\lambda$ will be the length of a graded piece.
\end{Example}

\begin{Lemma}\label{rank-lemma}
Let $B$ be a $\mathbb{Z}^{d+1}$-graded differential $R[y]$-module and
suppose we have the following commutative diagrams in $\mathcal{C}$:
\[
\xymatrix{
A\ar@{->>}_{\psi_A}[d]\ar^{\iota}[r] & B\ar^{\psi_B}[d] &&
	B'\ar@{->>}^{\varepsilon'}[r]\ar_{\phi_B}[d] &
	C'\ar@{^{(}->}^{\phi_C}[d] \\
A''\ar@{^{(}->}_{\iota''}[r] & B'' && B\ar_{\varepsilon}[r] & C\\
}\]
Then for each $\mathbf{m}\in\mathbb{Z}^{d+1}$ the following inequalities
hold:
\[
\lambda((\im\iota)_\mathbf{m})\geq\lambda((\im\iota'')_\mathbf{m})
\qquad\text{and}\qquad
\lambda((\im\varepsilon)_\mathbf{m})\geq\lambda(\im\varepsilon')_\mathbf{m}.
\]
Furthermore, if $\varepsilon\iota=0$ then
\[\lambda(B_\mathbf{m})\geq\lambda((\im\iota'')_\mathbf{m})+
\lambda((\im\varepsilon')_\mathbf{m}).\]
\end{Lemma}

\begin{proof}
For the first inequality, there is a surjection
$\im\iota\twoheadrightarrow \im\psi_B\iota$, so
\[\lambda((\im\iota)_\mathbf{m})\geq \lambda((\im\psi_B\iota)_\mathbf{m})=
\lambda((\im\iota''\psi_A)_\mathbf{m}).\]
Because $\psi_A$ is surjective there is also a surjection
$\im\iota''\psi_A\twoheadrightarrow\im\iota''$. This gives the
desired inequality,
$\lambda((\im\iota)_\mathbf{m})\geq\lambda((\im\iota'')_\mathbf{m})$.

For the second inequality, there is an inclusion
$\im\varepsilon'\hookrightarrow \im\varepsilon$ since $\phi_C$ is
injective. By superadditivity,
$\lambda((\im\varepsilon)_\mathbf{m})\geq\lambda((\im\varepsilon')_\mathbf{m})$.

For the final inequality, note that $\varepsilon\iota=0$ implies that
$\im\iota\subseteq\ker\varepsilon$. The exact sequence
\[\xymatrix{0\ar[r] & \ker\varepsilon\ar[r] & B\ar[r] &
\im\varepsilon\ar[r] & 0},\]
then implies $$\lambda(B_\mathbf{m})\geq
\lambda((\im\varepsilon)_\mathbf{m})+\lambda((\im\iota)_\mathbf{m})\geq
\lambda((\im\varepsilon')_\mathbf{m})+\lambda((\im\iota'')_\mathbf{m})$$ using
the first two inequalities.
\end{proof}

\begin{Lemma}\label{seq-lemma}
Let $D$ be a $\mathbb{Z}^{d+1}$-graded differential $R[y]$-module. Viewing
$R[y]\boxtimes_R D$ as a $R[y]$-module via the action $r(s\otimes
d)=(rs)\otimes d$, there is a sequence of $\mathbb{Z}^{d+1}$-graded
differential $R[y]$-modules
\[\xymatrix{0\ar[r] & (R[y]\boxtimes_R D)(-\deg y) \ar^-{\sigma}[r] &
R[y]\boxtimes_R D \ar^-{\varepsilon}[r] & D \ar[r] & 0},\]
with $\sigma(1\otimes d) = y\otimes d - 1\otimes yd$ and
$\varepsilon(a\otimes d)=ad$.
This sequence is exact and functorial in $D$. The map $\sigma$ is given by
multiplication by $y$ if and only if $yD=0$.
\end{Lemma}
\begin{proof}
It is straight-forward to check that $\sigma$ and $\varepsilon$ are
morphisms and that the sequence is exact and functorial.
Evidently $\sigma$ is multiplication by $y$ when $yD=0$. The exactness of
the sequence shows that the converse holds.
\end{proof}

The following definition and theorem are differential module versions
of Santoni's results for $R$-modules \cite{MR1028267}.

\begin{Definition}\label{single-var-high-low}
Let $A = k[x_1, \ldots, x_d]$.  A $\mathbb{Z}^{d+1}$-graded differential $A[y]$-module $D$ admits a
\emph{high-low decomposition} if there are non-zero
$\mathbb{Z}^{d+1}$-graded differential $A[y]$-modules $D_h$ and $D_\ell$
each annihilated by $y$, and there are morphisms of differential
$A[y]$-modules $\xymatrix{D_h\ar@{^{(}->}[r] & D}$ and
$\xymatrix{D\ar@{->>}[r] & D_\ell}$ (injective and surjective, respectively) that split in the category of
$\mathbb{Z}^d$-graded differential $A$-modules.

More generally if $R = k[x_1, \ldots, x_{n+1}]$ is a polynomial ring of dimension $n+1$, then we say that a $\mathbb{Z}^{n+1}$-graded differential $R$-module $D$ admits a \emph{high-low decomposition} if there is a variable $x_i$ ($1\leq i \leq n$) such that if $A = k[x_1, \ldots, \widehat{x_i}, \ldots, x_{n+1}]$ and $y$ is the indeterminate $x_i$ then $D$, viewed as a differential $A[y]$-module has a high-low decomposition as defined above.
\end{Definition}

\begin{Remark}
The difference between the two definitions above is minor, and is tantamount to a relabeling of the variables.  This is necessary for us to inductively apply high-low decompositions effectively in Lemma {\ref{bdd->high-low}}.  In the remaining results in Section 5, we will use only the first definition of high-low decompositions and remark that when we refer to $R[y]$ and \emph{high-low decompositions} then the identified additional variable is $y$. 
\end{Remark}

%\begin{Remark} \hl{The definition above follows that of {\cite{MR1028267}} but will frequently be applied in this paper in inductive cases where we regard $R = k[x_1, \ldots, \widehat{x_i}, \ldots, x_{d+1}]$ is and $x_i$ plays the role of $y$.   This is our approach in Lemma }\ref{bdd->high-low}.\end{Remark}

\begin{Theorem}\label{tor-thm}
Let $K$ be a $\mathbb{Z}^{d+1}$-graded $R[y]$-module such that $yK=0$, and
assume $\mathcal{C}$ is closed under $\Tor^{R[y]}(K,-)$. Let
$D\in\mathcal{C}$ be a $\mathbb{Z}^{d+1}$-graded differential module with
differential of degree $\mathbf{d}$ which admits a high-low decomposition.
Then for all $\mathbf{m}\in\mathbb{Z}^{d+1}$
\[
\lambda(\Tor^{R[y]}(K,D)_\mathbf{m})\geq
\lambda(\Tor^R(K,D_\ell)_{\mathbf{m}})+
\lambda(\Tor^R(K,D_h)_{\mathbf{m}+\mathbf{d}-\deg y}).
\]
\end{Theorem}
\begin{proof}
Applying the functoriality of Lemma~\ref{seq-lemma} to the high-low
decomposition
$\xymatrix{D_h\ar@{^{(}->}[r] & D}$ and $\xymatrix{D\ar@{->>}[r] & D_\ell}$
gives two exact commutative diagrams:
\[\xymatrix{
& 0\ar[d] & 0\ar[d] & 0\ar[d] \\
0\ar[r] & (R[y]\boxtimes_R D_h)(-\deg y) \ar^-{\sigma'}[r] \ar[d] &
R[y]\boxtimes_R D_h \ar^-{\varepsilon'}[r] \ar[d] & D_h \ar[r]\ar[d] & 0\\
0\ar[r] & (R[y]\boxtimes_R D)(-\deg y) \ar^-{\sigma}[r] &
R[y]\boxtimes_R D \ar^-{\varepsilon}[r] & D \ar[r] & 0\\}\]
and
\[\xymatrix{
0\ar[r] & (R[y]\boxtimes_R D)(-\deg y) \ar^-{\sigma}[r] \ar[d] &
R[y]\boxtimes_R D \ar^-{\varepsilon}[r] \ar[d] & D \ar[r] \ar[d] & 0\\
0\ar[r] & (R[y]\boxtimes_R D_\ell)(-\deg y) \ar^-{\sigma''}[r] \ar[d] &
R[y]\boxtimes_R D_\ell \ar^-{\varepsilon''}[r] \ar[d] & D_\ell \ar[r]\ar[d]
& 0\\
& 0 & 0 & 0\\}\]
In both diagrams the first two columns are split exact over $R[y]$ due to
the high-low decomposition. Because $D_h$ and $D_\ell$ are annihilated by
$y$, Lemma \ref{seq-lemma} implies that $\sigma'$ and $\sigma''$ are
multiplication by $y$. The $R[y]$-action on
$\Tor^{R[y]}(K,-)$ is via $K$ and $yK=0$, so after applying
$\Tor^{R[y]}(K,-)$ and Lemma \ref{change-of-rings} the maps $\sigma'$
and $\sigma''$ become zero, leaving
\begin{equation*}
\xymatrix@C-10.5pt{
& 0 \ar[d] && 0 \ar[d]\\
0\ar[r] & \Tor^R(K,D_h)\ar^{\varepsilon'}[r]\ar[d] &
\Tor^{R[y]}(K,D_h)\ar^-{\gamma'}[r]\ar[d] \ar@{}[dr] |{(\dag)}
& \Tor^R(K,D_h)(\mathbf{d}-\deg y)\ar[r]\ar[d] & 0\\
\dotso\ar^-{\sigma}[r] & \Tor^R(K,D)\ar^{\varepsilon}[r] &
\Tor^{R[y]}(K,D)\ar^-{\gamma}[r] & \Tor^R(K,D)(\mathbf{d}-\deg
y)\ar^-{\sigma(\mathbf{d})}[r]
& \dotso \\
}
\end{equation*}
and
\begin{equation*}
\xymatrix@C-10pt{
\dotso\ar^-{\sigma}[r] & \Tor^R(K,D)\ar^{\varepsilon}[r] \ar[d]
\ar@{}[dr] |{(\ddag)}
& \Tor^{R[y]}(K,D)\ar^-{\gamma}[r] \ar[d] & \Tor^R(K,D)(\mathbf{d}-\deg
y)\ar^-{\sigma(\mathbf{d})}[r]\ar[d] & \dotso\\
0\ar[r] & \Tor^R(K,D_\ell)\ar^{\varepsilon''}[r] \ar[d] &
\Tor^{R[y]}(K,D_\ell)\ar^-{\gamma''}[r] & \Tor^R(K,D_\ell)(\mathbf{d}-\deg
y)\ar[r]\ar[d] & 0\\
& 0 && 0\\
}
\end{equation*}
Lemma \ref{rank-lemma} on the commutative squares $(\dag)$ and
$(\ddag)$ completes the proof.
\end{proof}

\section{Lower bound on the Betti number}
\label{multi-grade-section}

In order to apply the results for high-low decompositions we need to
establish some results on the existence of high-low decompositions $D_h$
and $D_\ell$ with $H(D_h)\neq 0$ and $H(D_\ell)\neq 0$.

Recall that $\mathbf{m}_i$ denotes the $i$-th coordinate of a $d$-tuple
$\mathbf{m}\in\mathbb{Z}^d$.

\begin{Definition}
Let $D$ be a $\mathbb{Z}^d$-graded differential $R$-module and let $1\leq
i\leq d$. We say that $D$ is \emph{bounded in the $i$-th
direction} if there are $a,b\in\mathbb{Z}$ such that
$\mathbf{m}_i\not\in[a,b]$ implies $D_\mathbf{m}=0$.
\end{Definition}

\begin{Remark}\label{rmk:bdd<->finite-length}
When $D$ is finitely generated the condition that $D$ is bounded in the
$i$-th direction for all $i$ is equivalent to the condition that $\rank_k
D<\infty$.
\end{Remark}

\begin{Lemma}\label{bdd->high-low}
Let $D$ be a $\mathbb{Z}^d$-graded differential $R$-module with $H(D)\neq 0$.
Fix an index $1\leq i\leq d$ and suppose that $(\deg\delta^D)_i=0$.
If $H(D)$ is bounded in the $i$-th direction then there is a
$\mathbb{Z}^d$-graded differential module $D'$ that is quasi-isomorphic to
$D$ such that $D'$ has a high-low decomposition $D'_h$ and $D'_\ell$ with
$H(D'_h)$ and $H(D'_\ell)$ both non-zero.
\end{Lemma}

\begin{proof}
Let $a\in\mathbb{Z}$ be the largest integer such that $H(D)_\mathbf{m}=0$
whenever $\mathbf{m}_i<a$. Such an integer exists because $H(D)$ is non-zero
and bounded in the $i$-th direction. Set
\[E=\bigoplus_{\substack{\mathbf{m}\in\mathbb{Z}^d\\
a\leq\mathbf{m}_i}} D_\mathbf{m}.\]
This is an $R$-submodule. Since $(\deg\delta^D)_i=0$ it is closed under
$\delta^D$ as well. So $E$ is a differential submodule of $D$. By the
definition of $E$, we have
\[D/E \cong \bigoplus_{\substack{\mathbf{m}\in\mathbb{Z}^d\\ \mathbf{m}_i<a}}
D_\mathbf{m}.\]
We will use this isomorphism to identity with these differential modules (as a differential submodule of $D$)  in the calculations below.
Let $z$ be a cycle in $(D/E)_{\mathbf{m}}$. 
%removed by ref
%If $\mathbf{m}_i\geq a$ then
%$z=0$ as $(D/E)_{\mathbf{m}}=0$. 
Then $\mathbf{m}_i<a$ and 
$z\in (D/E)_\mathbf{m} \cong D_\mathbf{m}$ so there is a $z'\in D$ with
$\delta^D(z')=z$ as $H(D)_\mathbf{m}=0$ (here we use that $(\deg\delta^D)_i=0$).  So $\delta^{D/E}(z'+E)=z$.
Therefore $H(D/E)_\mathbf{m}=0$ for all $\mathbf{m}\in\mathbb{Z}^d$, and so
$H(D/E)=0$. From the short exact sequence $0\to E\to D \to D/E\to 0$ we conclude that $E\qis D$.

Let $b\in\mathbb{Z}$ be the smallest integer such that $H(E)_\mathbf{m}=0$
when $\mathbf{m}_i>b$. Again, such an integer exists because $H(E)\cong
H(D)$ is non-zero and bounded in the $i$-th direction. Set
\[E'=\bigoplus_{\substack{\mathbf{m}\in\mathbb{Z}^d\\
b+1\leq\mathbf{m}_i}} E_\mathbf{m}.\]
By an argument like the one above, $E'$ is a differential submodule of $E$ with $H(E')=0$ by the
definition of $b$. Set $D'=E/E'$.  As before, we conclude
that $H(E/E')\cong H(E)$ so that $D'=E/E'\qis E\qis D$.

By construction, $D'_\mathbf{m}=0$ for $\mathbf{m}_i<a$ (since $E$ is not supported in these degrees) and for
$\mathbf{m}_i>b$ (since everything in these degrees is in $E'$).  Also, by the definitions of $a$ and $b$, there are
$\mathbf{n},\mathbf{n}'\in\mathbb{Z}^d$ with $\mathbf{n}_i=a$ and
$\mathbf{n}'_i=b$ such that $H(D')_\mathbf{n}\neq 0$ and
$H(D')_{\mathbf{n}'}\neq
0$; hence $D'_\mathbf{n}\neq 0$ and $D'_{\mathbf{n}'}\neq 0$ as well.

Set
\[D'_\ell:= D'/\bigoplus_{\substack{\mathbf{m}\in\mathbb{Z}^d\\ \mathbf{m}_i > a}}
D'_\mathbf{m}
\cong \bigoplus_{\substack{\mathbf{m}\in\mathbb{Z}^d\\ \mathbf{m}_i=a}}
D'_\mathbf{m}
\qquad\text{and}\qquad
D'_h:=\bigoplus_{\substack{\mathbf{m}\in\mathbb{Z}^d\\ \mathbf{m}_i=b}}
D'_\mathbf{m} \subseteq D'. \]
Notice $D'_h$ is a submodule of $D'$ since multiplication by positive degree elements lands in a place where $D'$ is zero.  Then $D'_\ell$ and $D'_h$ are both non-zero and annihilated by $x_i$. The
two morphisms $\xymatrix{D'_h\ar@{^{(}->}[r] & D'}$ and
$\xymatrix{D'\ar@{->>}[r] & D'_\ell}$ split in the category of $\ZZ^{d-1}$-graded differential
modules because $(\deg\delta^D)_i=0$. So $D'_\ell$ and $D'_h$ form a
high-low decomposition. As noted above $H(D'_\ell)$ and $H(D'_h)$
are both non-zero, so $D'$ is the desired differential module.
\end{proof}

The proof of the following theorem uses Theorem~\ref{tor-thm}
inductively, after first using Lemma~\ref{bdd->high-low} to find a quasi-isomorphic
differential module with a suitable high-low decomposition.

Note that $H(D)$ is not required to be finitely generated in the following
theorem. If $H(D)$ is finitely generated then the hypothesis on $H(D)$ is
equivalent to $0<\rank_k H(D)<\infty$; see
Remark~\ref{rmk:bdd<->finite-length}.

\begin{Theorem}\label{total-Betti-bound}
If $D$ is a $\mathbb{Z}^d$-graded differential module
with $\deg\delta^D=\mathbf{0}$ and such that $H(D)\neq 0$ is bounded in the
$i$-th direction for all $i$, then
\[\beta^R(D)\geq 2^d.\]
\end{Theorem}

\begin{proof}
Use induction on $d$. For $d=0$, so that $R=k$, we have
\[\Tor^k(k,D)=H(k\boxtimes_k D)\cong H(D)\neq 0.\]
So $\beta^k(D)\geq 1$.

Now suppose $d>1$. Then $H(D)$ is bounded in the $d$-th direction by
assumption. By Proposition~\ref{diff-mod-tensor-exact} the Betti
number is preserved under quasi-isomorphisms, so
Lemma~\ref{bdd->high-low} allows us to assume that $D$ has a high-low
decomposition $D_h$ and $D_\ell$ with $H(D_h)\neq 0$ and $H(D_\ell)\neq 0$.
By definition of a high-low decomposition, $H(D_h)$ and $H(D_\ell)$ are
submodules of $H(D)$ over $k[x_1, \ldots, x_{d-1}]$ since the splitting happens in the category of
differential modules. In particular, $H(D_h)$ and $H(D_\ell)$ are bounded
in the $i$-th direction for all $i\in \mathbb{Z}^{d-1}$.

So the induction hypothesis applies to $D_h$ and $D_\ell$ thought of as $k[x_1,\ldots, x_{n-1}]$-modules. From
Theorem~\ref{tor-thm} we have:
\begin{align*}
\beta^R(D)&\geq \beta^{k[x_1,\dotsc,x_{d-1}]} (D_\ell)+
\beta^{k[x_1,\dotsc,x_{d-1}]}(D_h)\\
&\geq 2^{d-1}+2^{d-1}\\
&=2^d.\qedhere
\end{align*}
\end{proof}

\begin{Remark}
Example~\ref{scorpion} shows that Theorem~\ref{total-Betti-bound} cannot be
extended to differential modules $D$ with $\deg\delta^D>\mathbf{0}$.
\end{Remark}

Via Theorem~\ref{Betti-rank-thm} this result provides an
affirmative answer to Conjecture~\ref{hor-conj} when
$\deg\delta=\mathbf{0}$.

\begin{Corollary}
If $F$ is a finitely generated $\mathbb{Z}^d$-graded differential module
that is free as an $R$-module such that $\deg\delta^F=\mathbf{0}$ and such
that $H(F)$ has non-zero finite length then
\[\rank_R F\geq 2^d.\]
\end{Corollary}

\begin{proof}
By Theorem~\ref{non-pos-minimal-flag}, $F$ has a free flag. So
Theorem~\ref{Betti-rank-thm} implies that $\beta^R(F)\leq\rank_R F$.
Applying Theorem~\ref{total-Betti-bound} gives the desired inequality.
\end{proof}

\section*{Acknowledgments}
The authors gratefully acknowledge Srikanth Iyengar for his guidance, and
Luchezar Avramov and Lars Winther Christensen for helpful comments on this
paper. We also thank an anonymous referee for very helpful comments that
greatly improved this paper.

\bibliography{horrockref.bib}{}
\bibliographystyle{plain}
\end{document}